\documentclass[preprint,12pt]{elsarticle}

\usepackage{amssymb}
\usepackage{amsmath} 
\usepackage{mathrsfs}
\usepackage{amsmath,amssymb}
\usepackage{amsfonts}
\usepackage{latexsym,bm}
\usepackage{amsthm}
\usepackage{amssymb,amscd}
\usepackage{amsbsy}
\usepackage{fancyhdr,graphicx}
\usepackage[dvips]{psfrag}
\usepackage{indentfirst}
\usepackage{float}
\usepackage{graphics,color}
\usepackage{epsfig}
\usepackage{subfigure}
\usepackage{subcaption}
\usepackage{xcolor}
\usepackage{mdwlist}
\usepackage{epstopdf}
\usepackage{url}
\newtheorem {theorem} {Theorem} [section]

\newtheorem {lemma} {Lemma}[section]
\newtheorem {example} {Example}[section]
\newtheorem {remark} [theorem]{Remark}

\begin{document}

\begin{frontmatter}


 \tnotetext[t1]{The authors are partially supported NNSFC(Nos.12371501, 12071091, 12361034), Basic Science Research Project of Jiangsu higher education institutions, Ph. D. Fellowship of Xi’an Polytechnic University (No. 107020344), Science and Technology Plan Project of Guizhou Province(ZK[2022]G118).}
    \cortext[cor1]{Corresponding author}
\title{On the stability of singular Hopf bifurcation and its application}

\author[inst1]{Jun Li}
\affiliation[inst1]{organization={School of Mathmatics and Statistics},
            addressline={ Xidian University}, 
            city={Xi'an},
            postcode={710071}, 
            country={China}}
            \ead{lijun@xidian.edu.cn}
\author[inst2]{Shimin Li}
\affiliation[inst2]{organization={School of Mathmatics},
            addressline={Hangzhou Normal University}, 
            city={Hangzhou},
            postcode={311121}, 
            country={China}}
            \ead{lishimin@hznu.edu.cn}
\author[inst3]{Mingju Ma}
\affiliation[inst3]{organization={School of Science},
            addressline={Xi'an Polytechnic University}, 
            city={Xi'an},
            postcode={710048}, 
            country={China}}
            \ead{ mingjuma@163.com}
\author[inst4]{Kuilin Wu}
\ead{wukuilin@126.com}
\affiliation[inst4]{organization={School of Mathematics and Statistics},
            addressline={Guizhou University}, 
            city={Guiyang},
            postcode={550025}, 
            country={China}}

\begin{abstract}
Recently, research on the complex periodic behavior of multi-scale systems has become increasingly popular. Krupa et al. \cite{krupa2} provided a way to obtain relaxation oscillations in slow-fast  systems through singular Hopf bifurcations and canard explosion. The authors derived a $O(1)$ expression $A$ for the first Lyapunov coefficient (under the condition $A \neq 0$), and deduced the bifurcation curves of singular Hopf and canard explosions.

This paper employs Blow-up technique, normal form theory, and Lyapunov coefficient  formula to present higher-order approximate expressions for the first Lyapunov coefficient when $A=0$ for slow-fast systems. As an application, we investigate the bifurcation phenomena  of a predator-prey model with Allee effects. Utilizing the formulas obtained in this paper, we identify both supercritical and subcritical Hopf bifurcations that may occur simultaneously in the system. Numerical simulations validate the results. Finally, by normal form and slow divergence integral theory, we prove the cyclicity of the system  is 1. 
\end{abstract}




\end{frontmatter}


\section{Introduction}
In recent years, there has been a growing interest in the dynamic analysis of slow-fast  systems, which can be employed to describe dynamics for multiple time scales. There are  extensive applications  in various fields such as biochemistry \cite{41}, ecology \cite{li2023canard}, and engineering mechanics \cite{jardon2017model}. Geometric singular perturbation theory(GSPT), developed from the geometric theory pioneered by Fenichel \cite{F1}, serves as a fundamental framework for understanding such systems. Its core principle involves characterizing the original system as two distinct dynamical processes on fast and slow time scales. Through qualitative analysis of the layer and reduced systems, GSPT describes the dynamic behavior of the system when $\epsilon$ is sufficiently small.

In the early stages of GSPT, it was typically limited to situations where the critical manifold is normally hyperbolic. However, with the aid of Blow-up methods developed by Krupa \cite{krupa1,krupa2,krupa3}, De Maesschalck, Dumortier\cite{peter1,peter2,peter3,peter4,peter5} and others , it has become possible to derive normal forms near non-hyperbolic points on the critical manifold. This has led to the discovery of new bifurcation phenomena, such as canard explosions, singular Hopf bifurcations, and relaxation oscillations. These novel bifurcation phenomena have practical applications in population ecology \cite{41}, engineering mechanics \cite{7}, and other fields.

Krupa and Szmolyana\cite{krupa2} provide a detailed description for the process of canard explosion occurring at a fold point in the critical manifold. In this process, the generation of a singular Hopf bifurcation is of paramount importance. A Hopf bifurcation refers to the occurrence of  limit cycles due to the changes in stability of an equilibrium 
as a bifurcation parameters vary. If the limit cycle is stable (unstable), the bifurcation is termed supercritical (subcritical). In recent years, it has been observed in practical scenarios that, with changes in stability, multiple limit cycles may arise from a single equilibrium, even accompanied by the emergence of higher-codimensional degenerate singularity bifurcations \cite{huang1,huang2,huang3}. To determine the cyclicity of these periodic limit sets, it is necessary to compute Lyapunov coefficients. This analytical approach is widely applied in various systems.

In \cite{krupa2}, the authors provided a way to consider the normal form of the system and provided a $O(1)$ approximate expression $A$ for the first Lyapunov coefficient with respect to $\epsilon$. However, when $A=0$, we can not 
determine  whether Hopf bifurcation happens. It becomes essential to further compute higher-order expressions of the first Lyapunov coefficient with respect to $\epsilon$. This is crucial both for fundamental research and practical applications. In this paper, we aim to derive the formula for the first Lyapunov coefficient for planar differential systems under the multiscale framework. To achieve this, we consider a general slow-fast system:
\begin{equation}
    \begin{split}
        \frac{dx}{dt}&=f(x,y,\lambda,\epsilon),\\
        \frac{dy}{dt}&=\epsilon g(x,y,\lambda,\epsilon),
    \end{split}
    \label{d1}
\end{equation}
where, $f,g\in C^3$, $0 < \epsilon \ll 1$, and $\lambda$ is the bifurcation parameter. The system (\ref{d1}) is referred to as the fast system. For the fast system (\ref{d1}), after a timescale transformation $\tau = \epsilon t$, it is  transformed into the following slow system:
\begin{equation}
    \begin{split}
        \epsilon \frac{dx}{d\tau}&=f(x,y,\lambda,\epsilon),\\
        \frac{dy}{d\tau}&=g(x,y,\lambda,\epsilon).  
    \end{split}
    \label{d2}
\end{equation}
As $\epsilon \to 0$, the fast system (\ref{d1}) converges to the following layer system (\ref{d3}):
\begin{equation}
    \begin{split}
       \frac{dx}{dt}&=f(x,y,\lambda,\epsilon),\\
        \frac{dy}{dt}&=0,
    \end{split}
    \label{d3}
\end{equation}
and the slow system is reduced to
\begin{equation}
    \begin{split}
       &0=f(x,y,\lambda,\epsilon),\\
        &\frac{dy}{d\tau}=g(x,y,\lambda,\epsilon). 
    \end{split}
    \label{d4}
\end{equation}
The critical manifold can be defined by
$$
C=\left\{
(x,y)|f(x,y,\lambda,0)=0
\right\}.
$$
Without loss of generality, assume  that $O(0,0)$ is a local minimum point of the critical curve $C$, i.e.
\begin{equation*}
\begin{split}
    f(0,0,0,0)=0, ~~\frac{\partial f}{\partial x}(0,0,0,0)=0.  
\end{split}
\end{equation*}
We also need the following non-degeneracy assumption
\begin{equation*}
\begin{split}
    &\frac{\partial^2 f}{\partial x^2}(0,0,0,0)\ne 0, ~~\frac{\partial f}{\partial y}(0,0,0,0)\ne 0,\\
    &\frac{\partial g}{\partial x}(0,0,0,0)\ne 0, ~~\frac{\partial g}{\partial \lambda}(0,0,0,0)\ne 0.
\end{split}
\end{equation*}
Without loss of generality, suppose that
\begin{equation*}
\begin{split}
    &\frac{\partial^2 f}{\partial x^2}(0,0,0,0)> 0, ~~\frac{\partial f}{\partial y}(0,0,0,0)< 0,\\
    &\frac{\partial g}{\partial x}(0,0,0,0)> 0, ~~\frac{\partial g}{\partial \lambda}(0,0,0,0)< 0.
\end{split}
\end{equation*}
In this case, the system (\ref{d1}), when $(x, y, \lambda, \epsilon)$ is the neighborhood of $(0, 0, 0, 0)$, can be transformed into the following normal form through a series of diffeomorphisms (see  \cite{krupa2}) 
\begin{equation}
\begin{split}
    &\frac{dx}{dt}=-y h_1(x,y,\lambda,\epsilon)+x^2h_2(x,y,\lambda,\epsilon)+\epsilon h_3(x,y,\lambda,\epsilon),\\
    &\frac{dy}{dt}=\epsilon(xh_4(x,y,\lambda,\epsilon)-\lambda h_5(x,y,\lambda,\epsilon)+yh_6(x,y,\lambda,\epsilon),
\end{split}
\label{dnf}
\end{equation}
where
\begin{equation*}
    \begin{split}
&h_3(x,y,\lambda,\epsilon)=O(x,y,\lambda,\epsilon),\\
&h_i(x,y,\lambda,\epsilon)=1+O(x,y,\lambda,\epsilon),~~~i=1,2,4,5.
    \end{split}
\end{equation*}
Let
\begin{equation}
    \begin{split}
        &a_1=\frac{\partial h_3}{\partial x}(0,0,0,0),  a_2=\frac{\partial h_1}{\partial x}(0,0,0,0),   a_3=\frac{\partial h_2}{\partial x}(0,0,0,0),\\ 
       & a_4=\frac{\partial h_4}{\partial x}(0,0,0,0),   
        a_5=h_6(0,0,0,0),
    \end{split}
    \label{deffi4}
\end{equation}
and 
\begin{equation}
     A=-a_2+3 a_3-2 a_4-2 a_5.
\end{equation}
Denote $V$ by a small neighborhood of the origin. For singular Hopf bifurcations, there are the following classical results:
\begin{lemma}
\label{krupal}
    Suppose that the origin $O$ is a generic fold point for $ \lambda= 0$ with normal form (\ref{dnf}).  Then there exist $\epsilon_0>0,$ $ \lambda_0>0$ such that for each $0<\epsilon<\epsilon_0$, $|\lambda|<\lambda_0,$ equation (\ref{d1}) has precisely one equilibrium $p_e\in V$ which converges to the canard point as $(\epsilon,\lambda)\to 0.$ Moreover, there
exists a curve $\lambda_H(\sqrt{\epsilon})$ such that $p_e$ is stable for $\lambda<\lambda_H(\sqrt{\epsilon})$ and loses
stability through a Hopf bifurcation as $\lambda$ passes through $\lambda_H(\sqrt{\epsilon})$. The curve
$\lambda_H(\sqrt{\epsilon})$ has the expansion
$$
\lambda_H(\sqrt{\epsilon})=-\frac{a_1+a_5}{2}\epsilon+O(\epsilon^{3/2}).
$$
The Hopf bifurcation is non-degenerate if the constant $A$ defined in (\ref{deffi4})  is
nonzero. It is supercritical if $A<0$ and subcritical if $A>0.$
\end{lemma}
In practical applications, when $A=0$,  Lemma \ref{krupal} fail to provide relevant outcomes. In the subsequent sections of this paper, we will discuss the case $A=0$. To begin with, we present some relevant results(\cite{guk}) which can be applied to our results.
\begin{lemma} \label{DF}
    Consider the following planar differential system
\begin{equation}
    \begin{split}
        \frac{dx}{dt}=&-\beta_0 y+f(x,y),\\
        \frac{dy}{dt}=&\beta_0 x+g(x,y),
    \end{split}
    \label{deff3}
\end{equation}\label{dhp}
with $f(0,0)=g(0,0)=0,$ and $Df(0,0)=Dg(0,0)=0$, $\beta_0\ne 0,$ then the stability of the limit cycle is determined by the following first Lyapunov coefficient
\begin{equation}
    \begin{split}
        L_1&=\frac{1}{16}\{(f_{xxx}+f_{xyy}+g_{xxy}+g_{yyy})\\
        &~~~+\frac{1}{\beta_0}[f_{xy}(f_{xx}+f_{yy})-g_{xy}(g_{xx}+g_{yy})-f_{xx}g_{xx}+f_{yy}g_{yy}]\}|_{x=y=0}
    \end{split}
\end{equation}
A non-degenerate Hopf bifurcation is supercritical if the first Lyapunov coefficient
$L_1 < 0$, and subcritical if $L_1 > 0$.
\end{lemma}
To derive the expression for the first Lyapunov coefficient in the framework of slow-fast systems, we need the higher-order terms of $h_i$ in the normal form (\ref{dnf}). To ensure the uniqueness  of  system \eqref{dnf}, we need to rewrite system (\ref{dnf}) as follows:
\begin{equation}
    \begin{split}
    &\frac{dx}{dt}=-y h_1(x,y,\lambda)+x^2h_2(x,\lambda)+\epsilon h_3(x,y,\lambda,\epsilon),\\
    &\frac{dy}{dt}=\epsilon(xh_4(x,\epsilon)-\lambda h_5(x,y,\lambda,\epsilon)+yh_6(x,y,\epsilon),
\end{split}
\label{dnf2}
\end{equation}
where
\begin{equation*}
    \begin{split}
h_1(x,y,\lambda)&=1+\sum_{i+j=1}^{2}a_{ij}(\lambda)x^iy^j+O(|x,y|^3),\\
h_2(x,\lambda)&=1+b_{1,0}(\lambda )x+O(|x|^2),\\
h_3(x,y,\lambda,\epsilon)&=\sum_{i+j=1}^{3}c_{ij}(\lambda,\epsilon)x^iy^j+O(|x,y|^4),\\
h_4(x,\epsilon)&=1+\sum_{i+j=1}^{2}d_{i0}(\epsilon)x^i+O(|x|^3),\\
h_5(x,y,\lambda,\epsilon)&=1+\sum_{i+j=1}^{3}e_{ij}(\lambda,\epsilon)x^iy^j+O(|x,y|^3),\\
h_6(x,y,\epsilon)&=\sum_{i+j=0}^{2}f_{ij}(\epsilon)x^iy^j+O(|x,y|^3), ~~~~i,j\in N.\\
\end{split}
\end{equation*}
In the following sections, we will employ Lemma \ref{dhp} to further generalize the results in Lemma \ref{krupal}. 
\section{Main results and its proof}
\begin{theorem}\label{main}
 Suppose that the origin $O$ is a generic fold point for $ \lambda= 0$ with normal form (\ref{dnf2}).  Then there exist $\epsilon_0>0,$ $ \lambda_0>0$ such that for each $0<\epsilon<\epsilon_0$, $|\lambda|<\lambda_0,$ equation (\ref{d1}) has precisely one equilibrium $p_e\in V$ which converges to the canard point as $(\epsilon,\lambda)\to 0.$ Moreover, there
exists a curve $\lambda_1(\sqrt{\epsilon})$ such that $p_e$ is stable for $\lambda<\lambda_1(\sqrt{\epsilon})$ and loses
stability through a Hopf bifurcation as $\lambda$ passes through $\lambda_1(\sqrt{\epsilon})$. The curve
$\lambda_1(\sqrt{\epsilon})$ has the expansion
\begin{equation}
\lambda_1(\sqrt{\epsilon})=\epsilon^{1/2}\left(\rho_1 +\rho_3 \epsilon+O(\epsilon^{3/2})\right),
   \label{lambda*}
\end{equation}
where
$$ \rho _1=-\frac{1}{2} \left(c_{1,0}+f_{0,0}\right),~~~
        \rho_3=\frac{\rho _1}{8}(\rho_{31}+\rho _1\rho_{32}), 
$$
and 
\begin{equation}
    \begin{split}
       \rho_{31}&=a_{1,0} c_{1,0}+2 c_{1,0} f_{0,0}-2 c_{2,0}+e_{0,1}-f_{1,0},\\
       \rho_{32}&=a_{1,0}-3 b_{1,0}+2 \left(d_{1,0}-e_{1,0}+f_{0,0}\right).
    \end{split}
    \label{lambdaCoeff}
\end{equation}
The formula for the first Lyapunov coefficient of a singular Hopf bifurcation is:
\begin{equation}
        L_1(\sqrt{\epsilon})=\epsilon^{1/2}\left(\frac{1}{16}\omega_1 +\frac{1}{32} \omega_2 \epsilon+O(\epsilon^{3/2})\right),
        \label{LL}
\end{equation}
where
\begin{equation*}
\begin{split}
   &\omega_1= -a_{1,0}+3 b_{1,0}-2 (d_{1,0}+f_{0,0}),\end{split}
\end{equation*}
\begin{equation*}
    \begin{split}
   &\omega_2=6 a_{1,0} b_{1,0} c_{1,0}+6 a_{1,0} b_{1,0} f_{0,0}-4 a_{1,0} c_{1,0} d_{1,0}+a_{1,0} c_{1,0} e_{1,0}-4 a_{1,0} c_{1,0} f_{0,0}\\
   &~~~~-4 a_{1,0} c_{0,1}-2 a_{1,0}^2 c_{1,0}+2 a_{2,0} c_{1,0}-2 a_{1,0} c_{2,0}-6 a_{1,0} d_{1,0} f_{0,0}+a_{1,0} e_{1,0} f_{0,0}\\
   &~~~~-12 a_{1,0} f_{0,0}^2-4 a_{1,0}^2 f_{0,0}+6 a_{2,0} f_{0,0}+2 a_{0,1} \left(a_{1,0}+2 f_{0,0}\right)-2 a_{1,1}+2 f_{2,0}\\
   &~~~~+12 b_{1,0} c_{1,0} d_{1,0}-3 b_{1,0} c_{1,0} e_{1,0}+12 b_{1,0} c_{1,0} f_{0,0}+6 b_{1,0} c_{0,1}+12 b_{1,0} d_{1,0} f_{0,0}\\
   &~~~~-3 b_{1,0} e_{1,0} f_{0,0}+18 b_{1,0} f_{0,0}^2+4 c_{1,0} d_{1,0} e_{1,0}-8 c_{1,0} d_{1,0} f_{0,0}-8 c_{1,0} d_{1,0}^2-4 c_{0,1} d_{1,0}\\
   &~~~~-4 c_{2,0} d_{1,0}+6 c_{1,0} d_{2,0}+4 c_{1,0} e_{1,0} f_{0,0}-2 c_{1,0} e_{0,1}-2 c_{1,0} e_{2,0}-8 c_{0,1} f_{0,0}\\
   &~~~~-8 c_{2,0} f_{0,0}+2 c_{1,0} f_{1,0}+2 c_{1,1}+6 c_{3,0}+4 d_{1,0} e_{1,0} f_{0,0}-16 d_{1,0} f_{0,0}^2-8 d_{1,0}^2 f_{0,0}\\
   &~~~~+6 d_{2,0} f_{0,0}-2 d_{1,0} f_{1,0}+4 e_{1,0} f_{0,0}^2-2 e_{0,1} f_{0,0}-2 e_{2,0} f_{0,0}-8 f_{0,0}^3+4 f_{1,0} f_{0,0}.
\end{split}
\end{equation*}
The Hopf bifurcation is non-degenerate and 

(1)  it is supercritical if $\omega_1<0$ and subcritical if $\omega_1>0.$ 

(2) when $\omega_1=0$, the Hopf bifurcation is supercritical if $\omega_2|_{\omega_1=0}<0$ and subcritical if $\omega_2|_{\omega_1=0}>0.$
\end{theorem}
\begin{remark}
It is worth noting that, in the formula (\ref{LL}), $\omega_1$ is consistent with $A$ in theorem 3.1 and 3.2  \cite{krupa2}. However, when $A=0$, the direction of the Hopf bifurcation and the stability of the limit cycle can be determined by $\omega_2|_{\omega_1=0}$.
\end{remark}
\begin{proof}
We perform the following Blow-up transformation as follows:
$$
x=rx_1,y=r^2y_1,\lambda=r\lambda_1,\epsilon=r^2.
$$
Then system (\ref{dnf2})  can be transformed into:
\begin{equation}
    \begin{split}
        \frac{dx_1}{dt}&=\sum_{i+j=1}^3m_{i,j}x_1^i y_1^j+O(|x_1,y_1|^4),\\
        \frac{dy_1}{dt}&=n_{0,0}+\sum_{i+j=1}^3n_{i,j}x_1^i y_1^j+O(|x_1,y_1|^4),\\
    \end{split}
    \label{dnf3}
\end{equation}
where
\begin{equation*}
    \begin{split}
        m_{1,0}&=r c_{1,0}, m_{0,1}=-1+r^2 c_{0,1}, m_{2,0}=1+r^2 c_{2,0}, m_{1,1}=r \left(r^2 c_{1,1}-a_{1,0}\right),\\
         m_{0,2}&=r^2\left(r^2 c_{0,2}-a_{0,1}\right),m_{3,0}=r(b_{1,0}+r^2 c_{3,0}), m_{2,1}=r^2\left(r^2 c_{2,1}-a_{2,0}\right), \\
        m_{1,2}&=r^3\left(r^2 c_{1,2}-a_{1,1}\right), m_{0,3}=r^4\left(r^2 c_{0,3}-a_{0,2}\right), n_{0,0}=-\lambda _1, \\
        n_{1,0}&=1-\lambda _1 r e_{1,0}, n_{0,1}=r(f_{0,0}-\lambda _1 r e_{0,1}), n_{2,0}=r(d_{1,0}-\lambda _1 r e_{2,0}),\\
        n_{1,1}&=r^2(f_{1,0}-\lambda _1 r e_{1,1}),n_{0,2}=r^3(f_{0,1}-\lambda _1 r e_{0,2}), n_{3,0}=r^2(d_{2,0}-\lambda _1 r e_{3,0}), \\
        n_{2,1}&=r^3(f_{2,0}-\lambda _1 r e_{2,1}), n_{1,2}=r^4(f_{1,1}-\lambda _1 r e_{1,2}),n_{0,3}=r^5(f_{0,2}-\lambda _1 r e_{0,3}). 
    \end{split}
\end{equation*}
Assume that the equilibrium of the system \eqref{dnf3} is $p_r(x_1^*, y_1^*)$, which can be expressed as:
\begin{equation}
    \begin{split}
        x_1^*&=p_0+p_1r+p_2r^2+p_3r^3+O(r^4),\\ y_1^*&=q_0+q_1r+q_2r^2+q_3r^3+O(r^4),
    \end{split}
    \label{dnf4}
\end{equation}
where
\begin{equation*}
    \begin{split}
         &p_0=-\frac{n_{0,0}}{n_{1,0}}, p_1=-\frac{p_0^2 n_{2,0}+q_0 n_{0,1}}{n_{1,0}},  p_2=-\frac{p_0 \left(p_0^2 n_{3,0}+2 p_1 n_{2,0}+q_0 n_{1,1}\right)+q_1 n_{0,1}}{n_{1,0}},\\
         &p_3=-\frac{p_0 q_1 n_{1,1}+q_0 \left(p_0^2 n_{2,1}+p_1 n_{1,1}\right)+3 p_1 p_0^2 n_{3,0}+2 p_2 p_0 n_{2,0}+p_1^2 n_{2,0}+q_2 n_{0,1}+q_0^2 n_{0,2}}{n_{1,0}}, \\
         &q_0=-\frac{m_{2,0} n_{0,0}^2}{m_{0,1} n_{1,0}^2},\\
          &q_1=-\frac{p_0 \left(p_0^2 \left(m_{3,0} n_{1,0}-2 m_{2,0} n_{2,0}\right)+q_0 \left(m_{1,1} n_{1,0}-2 m_{2,0} n_{0,1}\right)+m_{1,0} n_{1,0}\right)}{m_{0,1} n_{1,0}},  \\
         &q_2=\frac{1}{m_{0,1} n_{1,0}}\times  \\
         &\left(p_0^2 \left(p_1 \left(4 m_{2,0} n_{2,0}-3 m_{3,0} n_{1,0}\right)+q_0 \left(2 m_{2,0} n_{1,1}-m_{2,1} n_{1,0}\right)\right)+p_0 q_1 \left(2 m_{2,0} n_{0,1}-m_{1,1} n_{1,0}\right)\right.\\
         &\left.-n_{1,0} \left(p_1 q_0 m_{1,1}+p_1 \left(p_1 m_{2,0}+m_{1,0}\right)+q_0^2 m_{0,2}\right)+2 p_0^4 m_{2,0} n_{3,0}\right),\\     
            \end{split}
        \end{equation*}
        \begin{equation*}
            \begin{split}
         &q_3=\frac{1}{m_{0,1} n_{1,0}}\times\\
         &\left(2 p_0^3 m_{2,0} \left(3 p_1 n_{3,0}+q_0 n_{2,1}\right)+p_0^2 \left(p_2 \left(4 m_{2,0} n_{2,0}-3 m_{3,0} n_{1,0}\right)+q_1 \left(2 m_{2,0} n_{1,1}-m_{2,1} n_{1,0}\right)\right)\right.\\
         &\left.+p_0 \left(2 p_1 q_0 \left(m_{2,0} n_{1,1}-m_{2,1} n_{1,0}\right)+p_1^2 \left(2 m_{2,0} n_{2,0}-3 m_{3,0} n_{1,0}\right)+q_2 \left(2 m_{2,0} n_{0,1}-m_{1,1} n_{1,0}\right)\right.\right.\\
         &\left.\left.+q_0^2 \left(2 m_{2,0} n_{0,2}-m_{1,2} n_{1,0}\right)\right)-n_{1,0} \left(p_1 q_1 m_{1,1}+q_0 \left(p_2 m_{1,1}+2 q_1 m_{0,2}\right)+p_2 \left(2 p_1 m_{2,0}+m_{1,0}\right)\right)\right).
    \end{split}
\end{equation*}
By transforming $x_2=x_1-x_1^*$ and $y_2=y_1-y_1^*$ to move  $p_r(x_1^*,y_1^*)$ to the origin,  system (\ref{dnf3}) can be transformed into:
\begin{equation}
    \begin{split}
        \frac{dx_2}{dt}&=\sum_{i+j=1}^3\bar{m}_{i,j}x_2^i y_2^j+O(|x_2,y_2|^4),\\
        \frac{dy_2}{dt}&=\sum_{i+j=1}^3\bar{n}_{i,j}x_2^i y_2^j+O(|x_2,y_2|^4),\\
    \end{split}
    \label{dnf5}
\end{equation}
where
\begin{equation*}
    \begin{split}
        \bar{m}_{1,0}&=2 p_0 m_{2,0}+\left(m_{1,0}+q_0 m_{1,1}+2 p_1 m_{2,0}+3 p_0^2 m_{3,0}\right)r\\
        &~+\left(q_1 m_{1,1}+2 p_2 m_{2,0}+2 p_0 q_0 m_{2,1}+6 p_0 p_1 m_{3,0}\right) r^2\\
        &~+\left(m_{1,2} q_0^2+q_2 m_{1,1}+2 p_3 m_{2,0}+2 \left(p_1 q_0+p_0 q_1\right) m_{2,1}\right.\\
        &~\left.+\left(3p_1^2+6 p_0 p_2\right) m_{3,0}\right) r^3+O\left(r^4\right), \\
    \bar{m}_{0,1}&=m_{0,1}+p_0m_{1,1}r+\left(p_0^2 m_{2,1}+p_1 m_{1,1}+2 q_0 m_{0,2}\right)r^2 \\
        &~+\left(2 p_0 q_0 m_{1,2}+p_2 m_{1,1}+2 p_0 p_1 m_{2,1}+2 q_1 m_{0,2}\right)r^3+O\left(r^4\right), \\
        \bar{m}_{2,0}&=m_{2,0}+3 p_0 m_{3,0} r+\left(q_0 m_{2,1}+3 p_1m_{3,0}\right)r^2+\left(q_1 m_{2,1}+3 p_2 m_{3,0}\right) r^3+O\left(r^4\right), \\
        \bar{m}_{1,1}&=m_{1,1}r+2 p_0 m_{2,1} r^2+\left(2 p_1 m_{2,1}+2 q_0 m_{1,2}\right) r^3+O\left(r^4\right),\\
        \bar{m}_{0,2}&=m_{0,2}r^2 +p_0m_{1,2} r^3+O\left(r^4\right), \\
        \bar{m}_{3,0}&=m_{3,0} r+O\left(r^4\right),  \bar{m}_{2,1}=m_{2,1} r^2+O\left(r^4\right),  \bar{m}_{1,2}= m_{1,2}r^3+O\left(r^4\right),\\  \bar{m}_{0,3}&=O\left(r^4\right), \\
        \bar{n}_{1,0}&=n_{1,0}+2 p_0n_{2,0}r+\left(3 p_0^2 n_{3,0}+2 p_1 n_{2,0}+q_0 n_{1,1}\right)r^2+\left(2 p_0 q_0 n_{2,1}+2 p_2 n_{2,0}\right.\\
        &~\left.+6 p_0 p_1 n_{3,0}+q_1 n_{1,1}\right)r^3  +O\left(r^4\right), \\ 
        \bar{n}_{0,1}&=n_{0,1} r+p_0 n_{1,1} r^2+\left(n_{2,1} p_0^2+2 q_0 n_{0,2}+p_1 n_{1,1}\right) r^3+O\left(r^4\right), \\
        \bar{n}_{2,0}&=n_{2,0}r +3 p_0n_{3,0}r^2+\left(3 p_1 n_{3,0}+q_0 n_{2,1}\right)r^3+O\left(r^4\right),\\
        \bar{n}_{1,1}&=n_{1,1}r^2+2 p_0n_{2,1}r^3+O\left(r^4\right), \\
        \bar{n}_{0,2}&=n_{0,2}r^3+O\left(r^4\right), 
        \bar{n}_{3,0}=n_{3,0}r^2+O\left(r^4\right),  \bar{n}_{2,1}= n_{2,1}r^3+O\left(r^4\right),\\
        \bar{n}_{1,2}&=O\left(r^4\right), \bar{n}_{0,3}=O(r)^5.
    \end{split}
\end{equation*}
Denote $M=-(\bar{m}_{1,0}+\bar{n}_{0,1}),N=\bar{m}_{1,0}\bar{n}_{0,1}-\bar{m}_{0,1}\bar{n}_{1,0}.$ \\
If $\bar{n}_{1,0}\ne 0$, we perform the following transformation:
$$
x_3=-\sqrt{2}\bar{n}_{1,0}x_2+\sqrt{2}\left(\bar{m}_{1,0}+\frac{M}{2}\right)y_2,y_3=\frac{\sqrt{2}}{2}\sqrt{4N-M^2}y_2.
$$
If $\bar{m}_{0,1} \ne 0$, we apply the transformation:
$$
x_3=\sqrt{2}\left(\bar{n}_{0,1}+\frac{M}{2}\right)x_2-\sqrt{2}\bar{m}_{0,1}y_2, y_3=\frac{\sqrt{2}}{2}\sqrt{4N-M^2}x_2.
$$
Thus,  system (\ref{dnf5}) can be reformulated in the following equivalent form:
\begin{equation}
    \begin{split}
        \frac{d x_3}{d t}&=\tilde{m}_{1,0}x_3- \tilde{m}_{0,1}y_3+\sum_{i+j=2}^{3}\tilde{m}_{i,j}x_{3}^i y_{3}^j+O(|x_3,y_3|^4),\\
        \frac{d y_3}{d t}&=\tilde{m}_{0,1}x_3+ \tilde{m}_{1,0}y_3+\sum_{i+j=2}^{3}\tilde{n}_{i,j}x_{3}^i y_{3}^j+O(|x_3,y_3|^4),
    \end{split}
    \label{dnf6}
\end{equation}
where
\begin{equation*}
    \begin{split}
        \tilde{m}_{1,0}&=\frac{1}{2} \left(\bar{m}_{1,0}+\bar{n}_{0,1}\right), \tilde{m}_{0,1}=-\frac{1}{2} \sqrt{2 \bar{m}_{1,0} \bar{n}_{0,1}-4 \bar{m}_{0,1} \bar{n}_{1,0}-\bar{m}_{1,0}^2-\bar{n}_{0,1}^2},\\
        \tilde{m}_{2,0}&=\frac{\bar{m}_{0,2} \left(\bar{n}_{0,1}-\bar{m}_{1,0}\right)-2 \bar{m}_{0,1} \bar{n}_{0,2}}{2 \sqrt{2} \bar{m}_{0,1}^2},\\
        \tilde{m}_{1,1}&=\frac{1}{\sqrt{2} \bar{m}_{0,1}^2 \sqrt{2 \bar{m}_{1,0} \bar{n}_{0,1}-4 \bar{m}_{0,1} \bar{n}_{1,0}-\bar{m}_{1,0}^2-\bar{n}_{0,1}^2}}\times\\
        &~~~\left(\bar{m}_{0,1} \left(-\bar{m}_{1,1} \bar{n}_{0,1}+\bar{m}_{1,0} \left(\bar{m}_{1,1}-2 \bar{n}_{0,2}\right)+2 \bar{m}_{0,1} \bar{n}_{1,1}+2 \bar{n}_{0,2} \bar{n}_{0,1}\right)\right.\\
        &~~~\left.-\bar{m}_{0,2} \left(\bar{m}_{1,0}-\bar{n}_{0,1}\right){}^2\right),\\
        \tilde{m}_{0,2}&=\frac{1}{2 \sqrt{2} \bar{m}_{0,1}^2 \left(-2 \bar{m}_{1,0} \bar{n}_{0,1}+4 \bar{m}_{0,1} \bar{n}_{1,0}+\bar{m}_{1,0}^2+\bar{n}_{0,1}^2\right)}\times\\
        &~~~\left(\bar{m}_{0,2} \left(\bar{m}_{1,0}-\bar{n}_{0,1}\right){}^3+2 \bar{m}_{0,1} \left(4 \bar{m}_{0,1}^2 \bar{n}_{2,0}+2 \bar{m}_{0,1} \bar{n}_{0,1} \left(\bar{n}_{1,1}-\bar{m}_{2,0}\right)\right.\right.\\
        &~~~\left.\left.+\bar{m}_{1,0}^2 \left(\bar{n}_{0,2}-\bar{m}_{1,1}\right)+\bar{n}_{0,1}^2 \left(\bar{n}_{0,2}-\bar{m}_{1,1}\right)+2 \bar{m}_{1,0} \left(\bar{n}_{0,1} \left(\bar{m}_{1,1}-\bar{n}_{0,2}\right)\right.\right.\right.\\
        &~~~\left.\left.\left.+\bar{m}_{0,1} \left(\bar{m}_{2,0}-\bar{n}_{1,1}\right)\right)\right)\right),\\
        \tilde{m}_{3,0}&=\frac{\bar{m}_{0,3} \left(\bar{m}_{1,0}-\bar{n}_{0,1}\right)+2 \bar{m}_{0,1} \bar{n}_{0,3}}{4 \bar{m}_{0,1}^3},\\
        \end{split}
        \end{equation*}
        \begin{equation*}
            \begin{split}
        \tilde{m}_{2,1}&=\frac{1}{4 \bar{m}_{0,1}^3 \sqrt{2 \bar{m}_{1,0} \bar{n}_{0,1}-4 \bar{m}_{0,1} \bar{n}_{1,0}-\bar{m}_{1,0}^2-\bar{n}_{0,1}^2}}\times\\
        &~~~\left(3 \bar{m}_{0,3} \left(\bar{m}_{1,0}-\bar{n}_{0,1}\right){}^2-2 \bar{m}_{0,1} \left(-\bar{m}_{1,2} \bar{n}_{0,1}+\bar{m}_{1,0} \left(\bar{m}_{1,2}-3 \bar{n}_{0,3}\right)\right.\right.\\
        &~~~\left.\left.+2 \bar{m}_{0,1} \bar{n}_{1,2}+3 \bar{n}_{0,3} \bar{n}_{0,1}\right)\right),\\
        \tilde{m}_{1,2}&=-\frac{1}{4 \bar{m}_{0,1}^3 \left(-2 \bar{m}_{1,0} \bar{n}_{0,1}+4 \bar{m}_{0,1} \bar{n}_{1,0}+\bar{m}_{1,0}^2+\bar{n}_{0,1}^2\right)}\times\\
        &~~\left(3 \bar{m}_{0,3} \left(\bar{m}_{1,0}-\bar{n}_{0,1}\right){}^3+2 \bar{m}_{0,1} \left(4 \bar{m}_{0,1}^2 \bar{n}_{2,1}-2 \bar{m}_{0,1} \bar{n}_{0,1} \left(\bar{m}_{2,1}-2 \bar{n}_{1,2}\right)\right.\right.\\
        &~~\left.\left.+\bar{m}_{1,0}^2 \left(3 \bar{n}_{0,3}-2 \bar{m}_{1,2}\right)+\bar{n}_{0,1}^2 \left(3 \bar{n}_{0,3}-2 \bar{m}_{1,2}\right)+2 \bar{m}_{1,0} \left(\bar{n}_{0,1} \left(2 \bar{m}_{1,2}-3 \bar{n}_{0,3}\right)\right.\right.\right.\\
        &~~\left.\left.\left.+\bar{m}_{0,1} \left(\bar{m}_{2,1}-2 \bar{n}_{1,2}\right)\right)\right)\right),\\
         \tilde{m}_{0,3}&=-\frac{1}{4 \bar{m}_{0,1}^3 \left(2 \bar{m}_{1,0} \bar{n}_{0,1}-4 \bar{m}_{0,1} \bar{n}_{1,0}-\bar{m}_{1,0}^2-\bar{n}_{0,1}^2\right){}^{3/2}}\times\\
        &~~\left(\left(2 \bar{m}_{0,1} \left(8 \bar{m}_{0,1}^3 \bar{n}_{3,0}+4 \bar{m}_{0,1}^2 \bar{n}_{0,1} \left(\bar{n}_{2,1}-\bar{m}_{3,0}\right)+2 \bar{m}_{0,1} \bar{n}_{0,1}^2 \left(\bar{n}_{1,2}-\bar{m}_{2,1}\right)\right.\right.\right.\\
        &~~\left.\left.\left.+\bar{m}_{1,0}^3 \left(\bar{m}_{1,2}-\bar{n}_{0,3}\right)+\bar{n}_{0,1}^3 \left(\bar{n}_{0,3}-\bar{m}_{1,2}\right)+\bar{m}_{1,0}^2 \left(3 \bar{n}_{0,1} \left(\bar{n}_{0,3}-\bar{m}_{1,2}\right)\right.\right.\right.\right.\\
        &~~\left.\left.\left.\left.-2 \bar{m}_{0,1} \left(\bar{m}_{2,1}-\bar{n}_{1,2}\right)\right)+\bar{m}_{1,0} \left(4 \bar{m}_{0,1}^2 \left(\bar{m}_{3,0}-\bar{n}_{2,1}\right)+4 \bar{m}_{0,1} \bar{n}_{0,1} \left(\bar{m}_{2,1}-\bar{n}_{1,2}\right)\right.\right.\right.\right.\\
        &~~\left.\left.\left.\left.+3 \bar{n}_{0,1}^2 \left(\bar{m}_{1,2}-\bar{n}_{0,3}\right)\right)\right)-\bar{m}_{0,3} \left(\bar{m}_{1,0}-\bar{n}_{0,1}\right){}^4\right.\right),\\
        \tilde{n}_{2,0}&=\frac{\bar{m}_{0,2} \sqrt{\bar{m}_{1,0} \bar{n}_{0,1}-2 \bar{m}_{0,1} \bar{n}_{1,0}-\frac{1}{2} \bar{m}_{1,0}^2-\frac{1}{2} \bar{n}_{0,1}^2}}{2 \bar{m}_{0,1}^2},\\
        \tilde{n}_{1,1}&=\frac{\bar{m}_{0,2} \left(\bar{m}_{1,0}-\bar{n}_{0,1}\right)-\bar{m}_{0,1} \bar{m}_{1,1}}{\sqrt{2} \bar{m}_{0,1}^2},\\
        \tilde{n}_{0,2}&=\frac{\bar{m}_{0,2} \left(\bar{m}_{1,0}-\bar{n}_{0,1}\right){}^2+2 \bar{m}_{0,1} \left(\bar{m}_{1,1} \bar{n}_{0,1}-\bar{m}_{1,0} \bar{m}_{1,1}+2 \bar{m}_{0,1} \bar{m}_{2,0}\right)}{2 \sqrt{2} \bar{m}_{0,1}^2 \sqrt{2 \bar{m}_{1,0} \bar{n}_{0,1}-4 \bar{m}_{0,1} \bar{n}_{1,0}-\bar{m}_{1,0}^2-\bar{n}_{0,1}^2}},\\
        \tilde{n}_{3,0}&=-\frac{\bar{m}_{0,3} \sqrt{2 \bar{m}_{1,0} \bar{n}_{0,1}-4 \bar{m}_{0,1} \bar{n}_{1,0}-\bar{m}_{1,0}^2-\bar{n}_{0,1}^2}}{4 \bar{m}_{0,1}^3},\\
        \tilde{n}_{2,1}&=\frac{2 \bar{m}_{0,1} \bar{m}_{1,2}-3 \bar{m}_{0,3} \left(\bar{m}_{1,0}-\bar{n}_{0,1}\right)}{4 \bar{m}_{0,1}^3},\\
        \tilde{n}_{1,2}&=\frac{-3 \bar{m}_{0,3} \left(\bar{m}_{1,0}-\bar{n}_{0,1}\right){}^2-4 \bar{m}_{0,1} \left(\bar{m}_{1,2} \bar{n}_{0,1}-\bar{m}_{1,0} \bar{m}_{1,2}+\bar{m}_{0,1} \bar{m}_{2,1}\right)}{4 \bar{m}_{0,1}^3 \sqrt{2 \bar{m}_{1,0} \bar{n}_{0,1}-4 \bar{m}_{0,1} \bar{n}_{1,0}-\bar{m}_{1,0}^2-\bar{n}_{0,1}^2}},\\
            \tilde{n}_{0,3}&=\frac{1}{4 \bar{m}_{0,1}^3 \left(-2 \bar{m}_{1,0} \bar{n}_{0,1}+4 \bar{m}_{0,1} \bar{n}_{1,0}+\bar{m}_{1,0}^2+\bar{n}_{0,1}^2\right)}\times
            \end{split}
        \end{equation*}
        \begin{equation*}
            \begin{split}
        &~~\left(\bar{m}_{0,3} \left(\bar{m}_{1,0}-\bar{n}_{0,1}\right){}^3-2 \bar{m}_{0,1} \left(2 \bar{m}_{2,1} \bar{m}_{0,1} \bar{n}_{0,1}+\bar{m}_{1,2} \bar{n}_{0,1}^2-2 \bar{m}_{1,0} \left(\bar{m}_{1,2} \bar{n}_{0,1}\right.\right.\right.\\
        &~~\left.\left.\left.+\bar{m}_{0,1} \bar{m}_{2,1}\right)+4 \bar{m}_{3,0} \bar{m}_{0,1}^2+\bar{m}_{1,0}^2 \bar{m}_{1,2}\right)\right).
    \end{split}
\end{equation*}
It is evident that the eigenvalues of Jacobian matrix at $(0,0)$  of system (\ref{dnf6})  is $\varrho = \tilde{m}_{1,0} \pm \tilde{m}_{0,1}i$. The necessary condition for the occurrence of a Hopf bifurcation in the system is $\tilde{m}_{1,0} = 0$, which is equivalent to:
\begin{equation*}
   \lambda_1=\lambda_1(r(\sqrt{\epsilon}))=\rho_1 r(\sqrt{\epsilon})+\rho_2 r^2(\sqrt{\epsilon})+\rho_3 r^3(\sqrt{\epsilon})+O(r^4(\sqrt{\epsilon})). 
\end{equation*}
Upon direct calculation, we find $\rho_2=0$, while $\rho_1$ and $\rho_3$ are determined by (\ref{lambdaCoeff}). Thus, when $\lambda_1=\lambda_1(r(\sqrt{\epsilon}))$, i.e., $\tilde{m}_{1,0}=0$, the system \eqref{dnf6} can be transformed into:
\begin{equation}
    \begin{split}
     \frac{d x_4}{d t}&=- \hat{m}_{0,1}y_4+\sum_{i+j=2}^{3}\hat{m}_{i,j}x_{3}^i y_{3}^j+O(|x_4,y_4|^4),\\
        \frac{d y_4}{d t}&=\hat{m}_{0,1}x_4+\sum_{i+j=2}^{3}\tilde{n}_{i,j}x_{3}^i y_{3}^j+O(|x_4,y_4|^4),
    \end{split}\label{dnf7}
\end{equation}
where
\begin{equation*}
    \begin{split}
        &\hat{m}_{0,1}=-\sqrt{-\bar{m}_{0,1} \bar{n}_{1,0}-\bar{m}_{1,0}^2},\hat{m}_{2,0}=-\frac{\bar{m}_{0,1} \bar{n}_{0,2}+\bar{m}_{0,2} \bar{m}_{1,0}}{\sqrt{2} \bar{m}_{0,1}^2},\\
        &\hat{m}_{1,1}=\frac{\bar{m}_{0,1} \left(\bar{m}_{1,0} \left(\bar{m}_{1,1}-2 \bar{n}_{0,2}\right)+\bar{m}_{0,1} \bar{n}_{1,1}\right)-2 \bar{m}_{0,2} \bar{m}_{1,0}^2}{\sqrt{2} \bar{m}_{0,1}^2 \sqrt{-\bar{m}_{0,1} \bar{n}_{1,0}-\bar{m}_{1,0}^2}},\\
        &\hat{m}_{0,2}=\frac{\bar{m}_{0,1} \left(\bar{m}_{0,1}^2 \bar{n}_{2,0}+\bar{m}_{1,0} \bar{m}_{0,1} \left(\bar{m}_{2,0}-\bar{n}_{1,1}\right)+\bar{m}_{1,0}^2 \left(\bar{n}_{0,2}-\bar{m}_{1,1}\right)\right)+\bar{m}_{0,2} \bar{m}_{1,0}^3}{\sqrt{2} \bar{m}_{0,1}^2 \left(\bar{m}_{0,1} \bar{n}_{1,0}+\bar{m}_{1,0}^2\right)},\\
        \hat{m}_{3,0}&=\frac{\bar{m}_{0,1} \bar{n}_{0,3}+\bar{m}_{0,3} \bar{m}_{1,0}}{2 \bar{m}_{0,1}^3},\hat{m}_{2,1}=\frac{3 \bar{m}_{0,3} \bar{m}_{1,0}^2-\bar{m}_{0,1} \left(\bar{m}_{1,0} \left(\bar{m}_{1,2}-3 \bar{n}_{0,3}\right)+\bar{m}_{0,1} \bar{n}_{1,2}\right)}{2 \bar{m}_{0,1}^3 \sqrt{-\bar{m}_{0,1} \bar{n}_{1,0}-\bar{m}_{1,0}^2}},\\
        &\hat{m}_{1,2}=\frac{-\bar{m}_{0,1} \left(\bar{m}_{0,1}^2 \bar{n}_{2,1}+\bar{m}_{1,0} \bar{m}_{0,1} \left(\bar{m}_{2,1}-2 \bar{n}_{1,2}\right)+\bar{m}_{1,0}^2 \left(3 \bar{n}_{0,3}-2 \bar{m}_{1,2}\right)\right)+3 \bar{m}_{0,3} \bar{m}_{1,0}^3}{2 \bar{m}_{0,1}^3 \left(\bar{m}_{0,1} \bar{n}_{1,0}+\bar{m}_{1,0}^2\right)},\\
        \hat{m}_{0,3}&=-\frac{1}{2 \bar{m}_{0,1}^3 \left(-\bar{m}_{0,1} \bar{n}_{1,0}-\bar{m}_{1,0}^2\right){}^{3/2}}\times\\
        &\left(\bar{m}_{0,1} \left(\bar{m}_{0,1}^3 \bar{n}_{3,0}+\bar{m}_{1,0} \bar{m}_{0,1}^2 \left(\bar{m}_{3,0}-\bar{n}_{2,1}\right)+\bar{m}_{1,0}^2 \bar{m}_{0,1} \left(\bar{n}_{1,2}-\bar{m}_{2,1}\right)+\bar{m}_{1,0}^3 \left(\bar{m}_{1,2}\right.\right.\right.\\
        &\left.\left.\left.\bar{n}_{0,3}\right)\right)-\bar{m}_{0,3} \bar{m}_{1,0}^4\right),\\
        \end{split}
        \end{equation*}
        \begin{equation*}
            \begin{split}
            \hat{n}_{2,0}&=\frac{\bar{m}_{0,2} \sqrt{-2 \bar{m}_{0,1} \bar{n}_{1,0}-2 \bar{m}_{1,0}^2}}{2 \bar{m}_{0,1}^2},\hat{n}_{1,1}=\frac{2 \bar{m}_{0,2} \bar{m}_{1,0}-\bar{m}_{0,1} \bar{m}_{1,1}}{\sqrt{2} \bar{m}_{0,1}^2},\\
        \hat{n}_{0,2}&=\frac{\bar{m}_{0,2} \bar{m}_{1,0}^2+\bar{m}_{0,1} \left(\bar{m}_{0,1} \bar{m}_{2,0}-\bar{m}_{1,0} \bar{m}_{1,1}\right)}{\sqrt{2} \bar{m}_{0,1}^2 \sqrt{-\bar{m}_{0,1} \bar{n}_{1,0}-\bar{m}_{1,0}^2}},\hat{n}_{3,0}=-\frac{\bar{m}_{0,3} \sqrt{-\bar{m}_{0,1} \bar{n}_{1,0}-\bar{m}_{1,0}^2}}{2 \bar{m}_{0,1}^3},\\
       \hat{n}_{2,1}&=\frac{\bar{m}_{0,1} \bar{m}_{1,2}-3 \bar{m}_{0,3} \bar{m}_{1,0}}{2 \bar{m}_{0,1}^3},\hat{n}_{1,2}=\frac{\bar{m}_{0,1} \left(2 \bar{m}_{1,0} \bar{m}_{1,2}-\bar{m}_{0,1} \bar{m}_{2,1}\right)-3 \bar{m}_{0,3} \bar{m}_{1,0}^2}{2 \bar{m}_{0,1}^3 \sqrt{-\bar{m}_{0,1} \bar{n}_{1,0}-\bar{m}_{1,0}^2}},\\
       \hat{n}_{0,3}&=\frac{\bar{m}_{0,3} \bar{m}_{1,0}^3-\bar{m}_{0,1} \left(\bar{m}_{3,0} \bar{m}_{0,1}^2-\bar{m}_{1,0} \bar{m}_{2,1} \bar{m}_{0,1}+\bar{m}_{1,0}^2 \bar{m}_{1,2}\right)}{2 \bar{m}_{0,1}^3 \left(\bar{m}_{0,1} \bar{n}_{1,0}+\bar{m}_{1,0}^2\right)}.
    \end{split}
\end{equation*}
System (\ref{dnf7}) is already in the required form (\ref{deff3}). By (\ref{lambda*}) and Lemma \ref{DF}, this completes the proof.
\end{proof}
\section{Applications}
In this section, we will utilize the results from Theorem \ref{main} to analyze the singular Hopf bifurcation problem in a predator-prey model with Allee effects.

In  \cite{biswas2023evolutionarily}, the authors analyzed a predator-prey model with the Allee effect in the prey's growth, which reveals that mating success at low densities is complicated due to difficulties in finding mates. In the model, this effect was described by the function $F(X) = \frac{rX}{\theta + X}$, where $X$ represents the population of the prey species, $r$ is the maximum per capita fertility rate, and $\theta$ represents the population density at which a species reaches half its maximum fertility, reflecting the strength of the Allee effect.

This model of the Allee effect was proposed by Ferdy \cite{ferdy2002allee}, which established a competition model including  the Allee effect in a patchy environment. The results showed that the Allee effect leaded to spatial segregation of species and maintains stability. In other words, populations with this effect can coexist in different spatial patches. Furthermore, this effect was applied to a predator-prey model with Holling Type II functional response \cite{zu2010impact}, . The authors found that as the Allee effect becomes stronger, the system may undergo subcritical Hopf bifurcations which leaded to unstable periodic oscillations. In other words, the system's equilibrium experienced a stability switch.
In  \cite{biswas2023evolutionarily}, this model of the Allee effect was applied to a predator-prey system in which the populations simultaneously have density-dependent terms. The model is as follows,
\begin{equation}
    \begin{split}
        \frac{dX}{dt}&=r \left(\frac{X}{\theta+X}\right)X-dX-k_1X^2-eXY,\\
        \frac{dY}{dt}&=\lambda  eXY-\delta Y-k_2Y^2,
    \end{split}
    \label{1}
\end{equation}
where $Y$ represents the population of predators, $d$ is the natural mortality rate, $k_1$ is the intra-species competition rate, and $e$ is the consumption rate. The function $F(x) = \frac{rX}{\theta + X}$ is used to describe the species fertility rate of the prey species $X$. $\lambda$, $\delta$ and $k_2$ represent the conversion rate, natural mortality rate and intra-species competition among predators, respectively.
The authors found that the system may exhibit various bifurcation phenomena. However, the existence conditions of limit cycles and other issues related to cyclicity for the system remain unresolved.

\begin{figure}
    \centering
\includegraphics[width=0.5\textwidth]{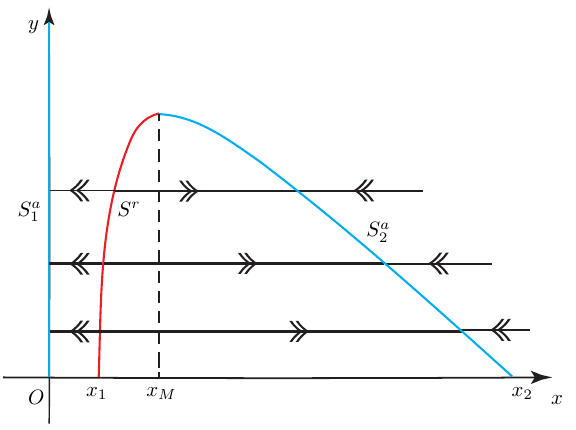}
    \caption{The critical manifold of system (\ref{2}).}
    \label{critical manifold}
\end{figure}

Assume that predator population has low conversion rate. And the average mortality,  the intra-species competition rate is also very low relative to prey population, i.e. $\lambda=\epsilon \Bar{\lambda},\delta=\epsilon \Bar{\delta}$ and  $k_2=\epsilon \Bar{k}_2$, $0<\epsilon \ll 1$.  Let $X = \frac{r}{k_1}x$, $Y = \frac{r}{e}y$, and $t = \frac{1}{r}\tau$, then system (\ref{1}) can be transformed into:
\begin{equation}
    \begin{split}
        \frac{dx}{d\tau}&=x\left(\frac{x}{m+x}-n-x-y\right):=f(x,y),\\
        \frac{dy}{d\tau}&=\epsilon y\left(\alpha x-\beta -\gamma y\right):=\epsilon g(x,y),
    \end{split}
    \label{2}
\end{equation}
where
$m=\frac{\theta k_1}{r}$, $n=\frac{d}{r}$, $\alpha=\frac{\Bar{\lambda}}{r}$, $\alpha=\frac{e\Bar{\lambda}}{k_1}$, $\beta=\frac{\Bar{\delta}}{r}$, $\gamma=\frac{\Bar{k_2}}{e}$. Note that $x=0$ and $y=0$ are two invariant lines of system (\ref{2}), and we state the following result without proof.
\begin{lemma}
The region $D=[0,1]\times[0,+\infty)$ is the forward invariant set of system (\ref{2}).
\end{lemma}
When $\epsilon=0$ in system (\ref{2}), we obtain the fast  subsystem:
\begin{equation}
    \begin{split}
        \frac{dx}{d\tau}&=x\left(\frac{x}{m+x}-n-x-y\right),\\
        \frac{dy}{d\tau}&=0.
    \end{split}
    \label{3}
\end{equation}
For system (\ref{2}), setting $\tau=s$ and $\epsilon=0$, we have the slow subsystem.
\begin{equation}
    \begin{split}
        \frac{dx}{ds}&=x\left(\frac{x}{m+x}-n-x-y\right),\\
        0&= y\left(\alpha x-\beta -\gamma y\right).
    \end{split}
    \label{4}
\end{equation}
Denote $F(x)=\frac{x}{m+x}-n-x$. By direct calculation, we can obtain the function $F(x)$ has a fold point $M(x_M, y_M)$ in the first quadrant when
\begin{equation}
\begin{split}
0 < m < \left(1 - \sqrt{n}\right)^2, \quad 0 < n < 1,
\end{split}
\label{5}
\end{equation}where $x_M = \sqrt{m} - m$ and $y_M = 1 - n + m - 2\sqrt{m}$.
Furthermore, denote:
\begin{equation*}
    \begin{split}
      S^a&=S_1^a\cup S_2^a,\\  S^r&=\left\{\left(x,F(x)\right)\mid \frac{1-m-n-\sqrt{\triangle}}{2}<x<\sqrt{m}-m\right\},  
    \end{split}
\end{equation*}
where
\begin{equation*}
    \begin{split}
    S_1^a&=\left\{\left(0,y\right)\mid y\ge 0\right\},\\
    S_2^a&=\left\{\left(x,F(x)\right)\mid \sqrt{m}-m<x<\frac{1-m-n+\sqrt{\triangle}}{2}\right\}.
\end{split}
\end{equation*}
The critical manifold can be divided into the normal attracting part $S^a$ and the normal repelling part $S^r,$   see Fig.\ref{critical manifold}.
\section{Equilibria of system (\ref{2})}
\begin{figure}
\centering
\begin{minipage}[t]{0.48\textwidth}
\centering
\includegraphics[width=6cm]{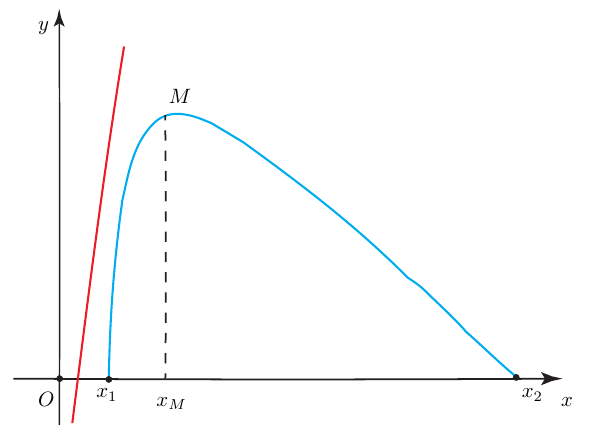}
\end{minipage}
\begin{minipage}[t]{0.48\textwidth}
\centering
\includegraphics[width=6cm]{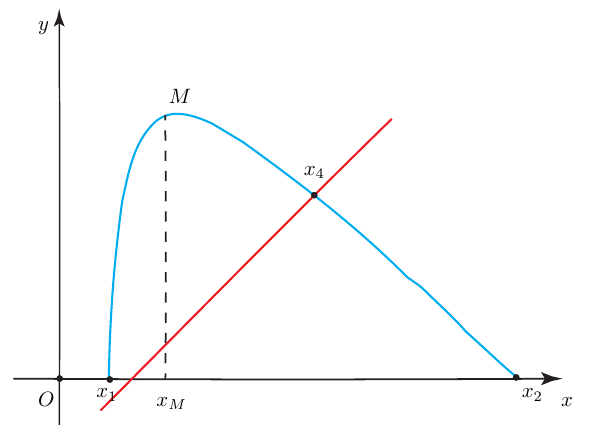}
\end{minipage}
\begin{minipage}[t]{0.48\textwidth}
\centering
\includegraphics[width=6cm]{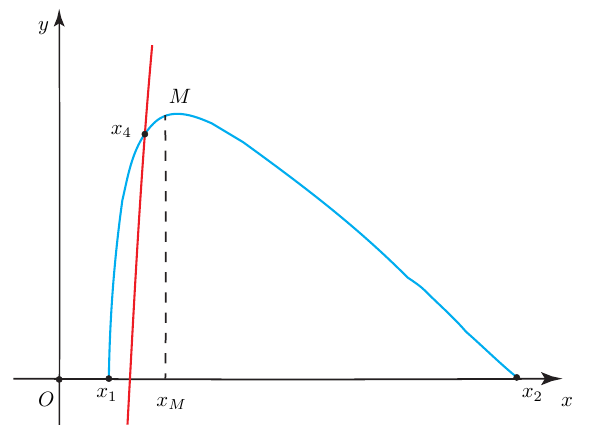}
\end{minipage}
\begin{minipage}[t]{0.48\textwidth}
\centering
\includegraphics[width=6cm]{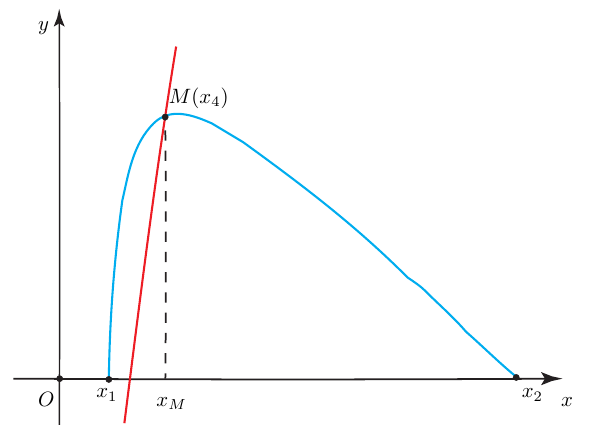}
\end{minipage}
\begin{minipage}[t]{0.48\textwidth}
\centering
\includegraphics[width=6cm]{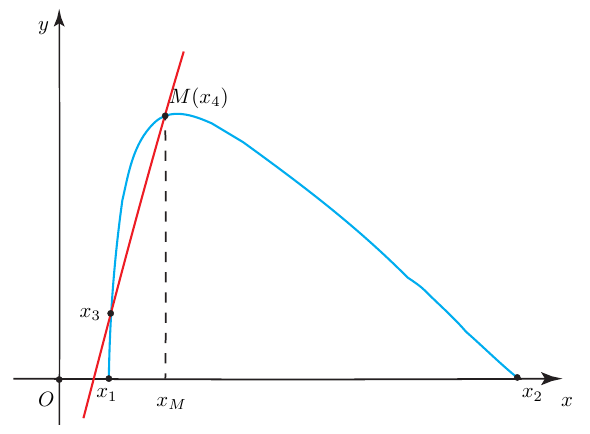}
\end{minipage}
\begin{minipage}[t]{0.48\textwidth}
\centering
\includegraphics[width=6cm]{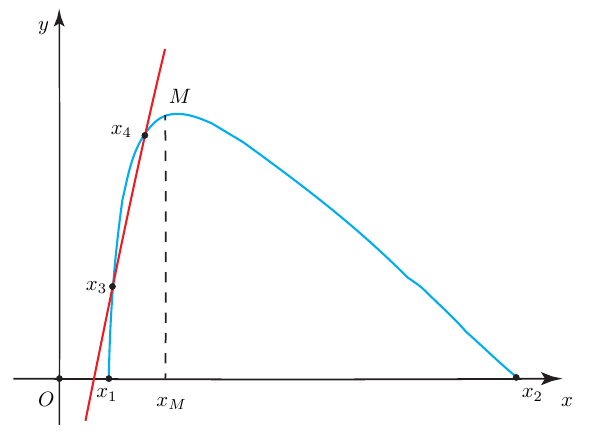}
\end{minipage}
\begin{minipage}[t]{0.48\textwidth}
\centering
\includegraphics[width=6cm]{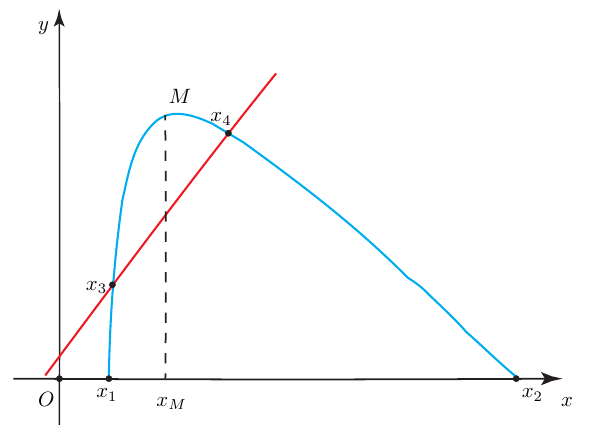}
\end{minipage}
\begin{minipage}[t]{0.48\textwidth}
\centering
\includegraphics[width=6cm]{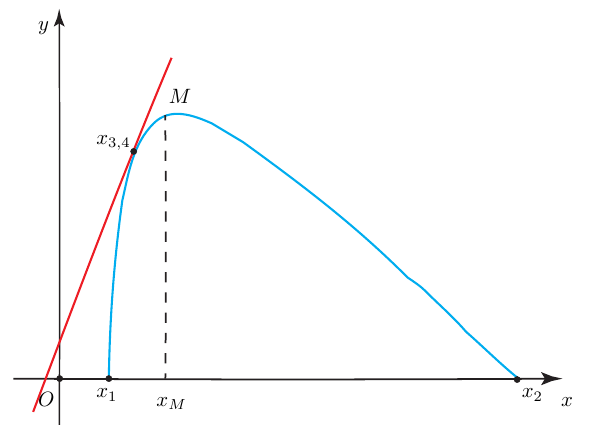}
\end{minipage}
\caption{The existence of equilibria $E_0, E_{1,2}, E_{3,4}$.}
\label{EE}
\end{figure}

System (\ref{2}) has an equilibrium  at $E_0(0,0)$. As for the boundary equilibria $E_{1,2}(x_{1,2},0)$, where $x_{1,2}$ satisfy the equation:
\begin{equation}
    x^2+(m+n-1)x+m n=0.
    \label{E12}
\end{equation}
If  $\bigtriangleup_1:=(1-m-n)^2-4m>0,$ and $1-m-n>0$, then $E_{1,2}$ exist (See Fig (\ref{EE})), and $x_{1,2}$ satisfy
\begin{equation*}
  x_{1,2}=\frac{1-m-n\mp \sqrt{\bigtriangleup_1}}{2}.
\end{equation*}
If $\bigtriangleup_1=0,$ $1-m-n>0$,$E_{1,2}$ will collide and become the unique boundary equilibrium $E^0_{1,2}$.
There are two  positive equilibria $E_{3,4}(x_{3,4},y_{3,4}),$ where
$y_{3,4}=\frac{\alpha x_{3,4}-\beta}{\gamma}$,
and $x_{3,4}=\frac{-b\mp\sqrt{\bigtriangleup_2}}{2}$ are the roots of the equation
\begin{equation}
x^2+bx+c=0,
\label{E34}
\end{equation}
where
$b=\frac{\gamma(m+n-1)+m \alpha-\beta}{\alpha+\gamma},c=\frac{m(\gamma n-\beta)}{\alpha+\gamma}.$
The equilibrium $E_{3,4}$ exist if and only if $\bigtriangleup_2=\frac{1}{\gamma+\alpha}\left( \gamma(m+n-1)+m\alpha-\beta\right)^2-4m(\alpha+\gamma)(\gamma n-\beta)>0$, as shown in Fig.(\ref{EE}).
\begin{lemma}
For system (\ref{2}), the boundary equilibrium $E_0(0,0)$ is a stable node. If $\bigtriangleup_1>0$ and $1-m-n>0$, then $E_1(x_1,0)$ and $E_2(x_2,0)$ exist, where $E_1$ and $E_2$ are saddles. If $\bigtriangleup_1>0$, $1-m-n>0$ and $\bigtriangleup_2>0$, then $E_3(x_3,y_3)$ and $E_4(x_4,y_4)$ also exist. In this case, $E_3$ is a saddle and $E_4$ is a non-saddle point.
\end{lemma}
\subsection{Singular Hopf bifurcation and Canard explosion}
\begin{figure}[ht]
\centering
\begin{minipage}[t]{0.48\textwidth}
\centering
\includegraphics[width=6cm]{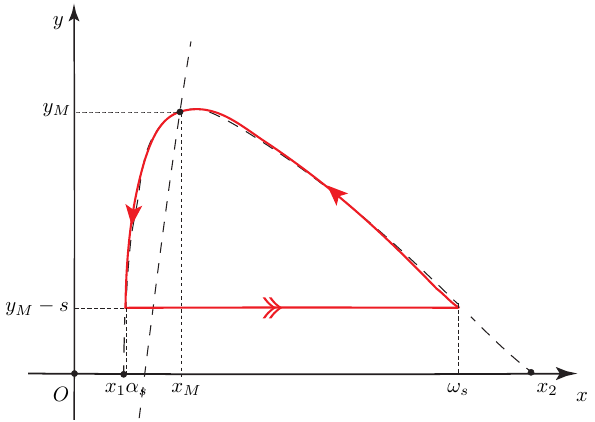}
\put(-80,80){$\Gamma(s)$}
\end{minipage}
\begin{minipage}[t]{0.48\textwidth}
\centering
\includegraphics[width=5cm]{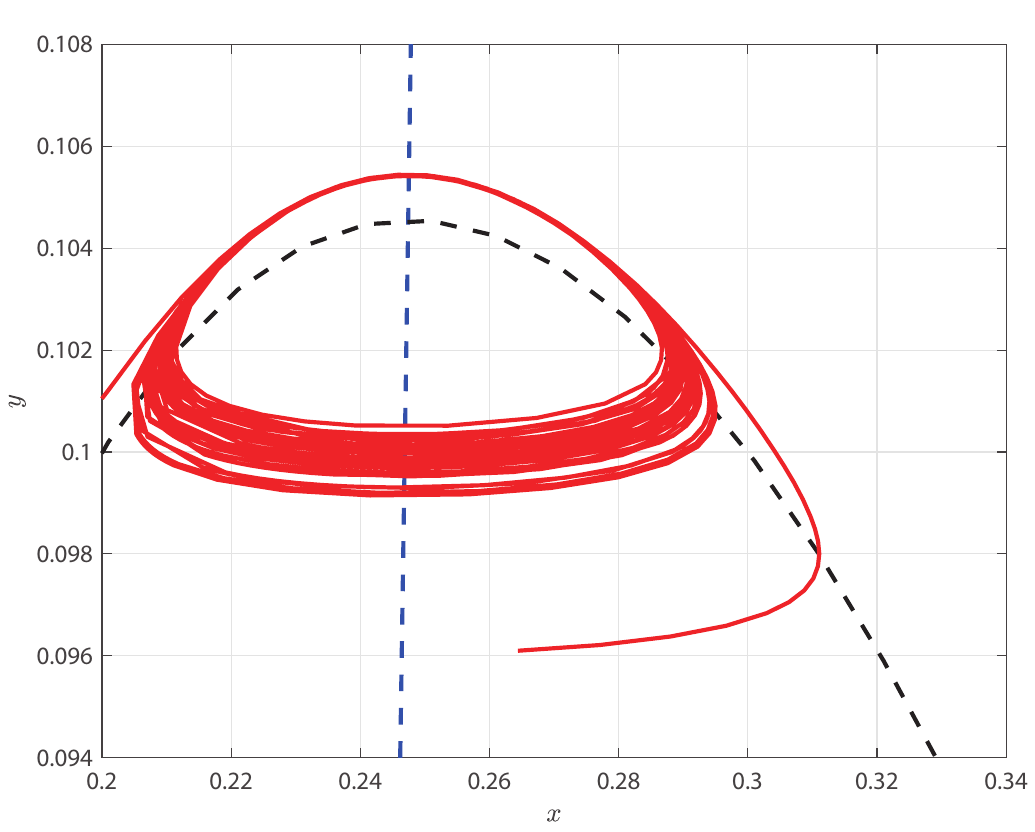}
\end{minipage}
\caption{If $E_4$ coincides with the critical point $M$, the system may have canard
explosion.}
\label{limit cycle}
\end{figure}

When the fold point $M(x_M, y_M)$ coincides with the positive equilibrium  $E_2$, i.e., $\beta + \gamma y_M - \alpha x_M = 0$, system  \eqref{2} may undergo a Hopf bifurcation (see Fig. (\ref{limit cycle})). To discuss the singular  Hopf bifurcation, we need to derive the normal form provided by Krupa et al.\cite{krupa2} and De Maesschalck et. al \cite{peter1,peter2,peter3,peter4,peter5}.
First, translate the point $M$ to the origin by  $\Bar{x}=x-x_M, \Bar{y}=y-y_M$, and make the transformations $\Bar{x}=\frac{\sqrt{\alpha x_M y_M}}{\sqrt{m}-1}X,$ $\Bar{y}=\frac{\alpha y_M}{\sqrt{m}-1}Y$, $\tau=\frac{1}{\sqrt{\alpha x_M y_M}}t$. Redefine $x=\bar{x},y=\bar{y}$ and let $\beta^*=\alpha  x_M-\gamma  y_M$, $\lambda = \beta -\beta _*$ system (\ref{2}) can be transformed into:
\begin{equation*}
    \begin{split}
        \frac{dx}{dt}&=-y h_1\left(x,y,\lambda,  \right)+x^2 h_2\left(x,\lambda,  \right)+\epsilon h_3\left(x,y,\lambda, \epsilon \right),\\
        \frac{dy}{dt}&=\epsilon\left(x h_4\left(x, \epsilon \right)-\lambda h_5\left(x,y,\lambda, \epsilon \right)+y h_6\left(x,y,\lambda  \right)\right),\\
    \end{split}
\end{equation*}
where 
\begin{equation*}
    \begin{split}
       h_1\left(x,y,\lambda, \epsilon \right)&=1+\frac{\alpha y_M}{\sqrt{\alpha x_M y_M}(\sqrt{m}-1)}x+O(|x,y|^4),\\
        h_2\left(x,y,\lambda, \epsilon \right)&=1-\frac{\sqrt{\alpha x_M y_M}}{(\sqrt{m}-1)^2}x+O(|x|^2),\\
        h_3\left(x,y,\lambda, \epsilon \right)&=0,\\
        h_4\left(x,y,\lambda, \epsilon \right)&=1++O(|x|^3),\\
         h_5\left(x,y,\lambda, \epsilon \right)&=1+\frac{\alpha }{\sqrt{m}-1}y+O(|x,y|^4),\\
         h_6\left(x,y,\lambda, \epsilon \right)&=-\frac{\gamma  y_M}{\sqrt{\alpha  x_M y_M}}+\frac{\alpha  x}{\sqrt{m}-1} +\frac{\gamma  \sqrt{\alpha  x_M y_M}}{(1-\sqrt{m}) x_M}y+O(|x,y|^3).\\
    \end{split}
\end{equation*}
If we choose $\lambda$ as the bifurcation parameter, then
\begin{equation*}
\begin{split}
\lambda= \frac{\left(\beta -\beta ^*\right) \left(1-\sqrt{m}\right)}{\alpha  \sqrt{\alpha  x_M y_M}} =0\Leftrightarrow \beta=\beta^*. 
    \end{split}
\end{equation*}

Next, we demonstrate the existence of the Hopf bifurcation. By directly computing, we can get
\begin{equation*}
    \begin{split}
        a_1&=\frac{\partial h_3}{\partial x}(0,0,0,0)=0,\\
        a_2&=\frac{\partial h_1}{\partial x}(0,0,0,0)=\frac{\alpha y_M}{\sqrt{\alpha x_M y_M}\sqrt{m}-1},\\
        a_3&=\frac{\partial h_2}{\partial x}(0,0,0,0)=-\frac{\sqrt{\alpha x_M y_M}}{(\sqrt{m}-1)^2},\\
        \end{split}
\end{equation*}
\begin{equation*}
    \begin{split}
        a_4&=\frac{\partial h_4}{\partial x}(0,0,0,0)=0,\\
        a_5&=h_6(0,0,0,0)=\frac{\alpha x_M-\beta-2\gamma y_M}{\sqrt{\alpha x_M y_M}},\\
        A&=-a_2+3a_3-2a_4-2a_5=\frac{\Psi(m)}{\sqrt{\alpha x_M y_M}(1-\sqrt{m})}y_M,
    \end{split}
\end{equation*}
where $\Psi(m)=2 \gamma(1-\sqrt{m})+\alpha-3\alpha \sqrt{m}.$ Note that $\frac{d\Psi (m)}{dm}=-\frac{3 \alpha +2 \gamma }{2 \sqrt{m}}<0$  and (\ref{5}), there are two cases as follows,

Case 1.  when $\frac{\alpha+2 \gamma}{3 \alpha+2 \gamma}\le 1-\sqrt{n}$, i.e. $n\le\left(\frac{2\alpha}{3 \alpha+2\gamma}\right)^2$, the following statements hold.  

a). if $0<m<\left(\frac{\alpha+2 \gamma}{3 \alpha+2 \gamma}\right)^2$, then $A>0.$

b). if $\left(\frac{\alpha+2 \gamma}{3 \alpha+2 \gamma}\right)^2<m<(1-\sqrt{n})^2$,then  $A<0.$

c). if $m=\left(\frac{\alpha+2 \gamma}{3 \alpha+2 \gamma}\right)^2$,then $A=0.$

Case 2. If $\frac{\alpha+2 \gamma}{3 \alpha+2 \gamma}>1-\sqrt{n}$, i.e. $n>\left(\frac{2\alpha}{3 \alpha+2\gamma}\right)^2,$  then $A>0.$

As $A=0$, it is easy to verify that the assumptions of Theorem \ref{main} is established. In this case, $m=\left(\frac{\alpha+2 \gamma}{3 \alpha+2 \gamma}\right)^2$, 
and 
\begin{equation*}
    \begin{split}\omega_2|_{\omega_1=0}&=\frac{\gamma  (3 \alpha +2 \gamma )^2 \left(9 \sqrt{2} (3 \alpha +2 \gamma ) \sqrt{\alpha +2 \gamma } y_M^{3/2}+8 \left(\sqrt{2} \alpha  \sqrt{(\alpha +2 \gamma ) y_M}+1\right)\right)}{8 \alpha ^2 (\alpha +2 \gamma )}\\
&>0.
\end{split}
\end{equation*}
According to Theorem \ref{main}, in this case, an unstable limit cycle can be bifurcated from  the fold point $M(x_M,y_M)$ of system \eqref{2}, and we have the following two results.
\begin{theorem}
For system (\ref{2}) with  $\gamma=\frac{\beta +\alpha  m-\alpha  \sqrt{m}}{-m+2 \sqrt{m}+n-1}$, $\bigtriangleup_1>0,$ $1-m-n>0$,  there is an equilibrium $E_4(x_4,y_4)$ near the fold point $M(x_M,y_M)$, which satisfies  $E_4\to M $ when $(\lambda,\epsilon)\to(0,0)$. Moreover, there is a singular Hopf bifurcation curve $\lambda_h(\sqrt{\epsilon})$, s.t.  the equilibrium $E_4$ is stable when $\lambda<\lambda_h(\sqrt{\epsilon})$.  When $\lambda$ passes through $\lambda_h(\sqrt{\epsilon})$, the system will undergo a Hopf bifurcation, where

\begin{equation*}
\lambda_h(\sqrt{\epsilon})=-\frac{a_1+a_5}{2}\epsilon+O(\epsilon^{3/2})=-\frac{\alpha x_M-\beta-2\gamma y_M}{2\sqrt{\alpha x_M y_M}}\epsilon+O(\epsilon^{3/2}).    
\end{equation*}
\end{theorem}

Applying the normal form (3.1), (3.15) and (3.16) in \cite{krupa2}, the singular Hopf bifurcation curve $\lambda_h(\epsilon)$ is obtained. Applying Theorem 3.2 in \cite{krupa2}, we have the following  Hopf bifurcation and canard explosion.

\begin{theorem}
For system (\ref{2}) with $\gamma=\frac{\beta +\alpha  m-\alpha  \sqrt{m}}{-m+2 \sqrt{m}+n-1}$, $\bigtriangleup_1>0,$ and $1-m-n>0$,   a $\epsilon$-family of canard cycles without head $\Gamma(\epsilon,s)$ bifurcate from the limit periodic set $\Gamma(s)$, where $s\in (0, y_M)$ if $y_3<0$, and $ s\in (0, y_M-y_3)$ if $y_3 \ge 0$. $\lambda=\lambda(s,\sqrt{\epsilon})$, $0<\epsilon\ll 1$. Moreover,  $\lambda(s,\sqrt{\epsilon})$ satisfies
\begin{equation*}
    |\lambda(s,\sqrt{\epsilon})-\lambda_c(\sqrt{\epsilon})|\le e^{-1/\epsilon^K},
\end{equation*}
where $K>0$ is a constant, and
\begin{equation}
\begin{split}
\lambda_c(\sqrt{\epsilon})&=-\left(\frac{a_1+a_5}{2}+\frac{A}{8}\right)\epsilon+O(\epsilon^{3/2})\\
    &=\frac{4 \left(\sqrt{m}-1\right) \left(\beta -\alpha  x_M\right)-y_M \left(\alpha +3 \alpha  \sqrt{m}-6 \gamma  \left(\sqrt{m}-1\right)\right)}{8 \left(\sqrt{m}-1\right) \sqrt{\alpha  x_M y_M}}\epsilon+O(\epsilon^{3/2}).
\end{split}
\label{hopf}
\end{equation}
\end{theorem}
\begin{example}
    Let $m=0.3,n=0.1,\gamma =0.1, \beta =0.2,\alpha =0.849561$, $\epsilon=0.0099$ and the initial value $(0.2644, 0.0961)$, then there is a small stable Hopf cycle around $E_4$, see Fig. \ref{limit cycle}(b).
\end{example}
\begin{figure}[ht]
    \centering
\includegraphics[width=0.5\textwidth]{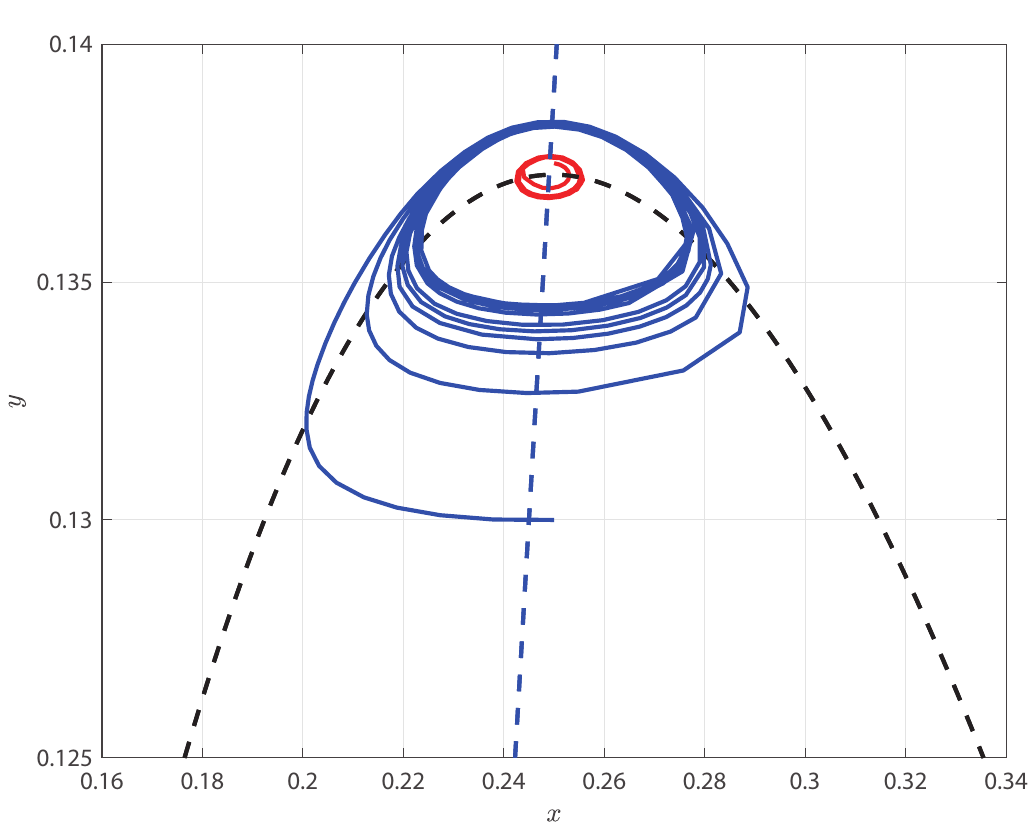}
    \caption{When $\alpha=0.8, \gamma=0.4424, n=0.1, \beta=0.138485, m=0.263075, \epsilon=0.01,$ and $A=0$. Apply a time reversal transformation $t\to -t,$ and take  initial value $(0.25, 0.1375)$ and $(0.25, 0.13),$ according to Poincaré-Bendixson theorem, system (\ref{2}) exhibits an unstable limit cycle.}
    \label{unstable}
\end{figure}
\begin{example}
 We apply a time reversal transformation to system (\ref{2}), $t\to -t,$ and let $\alpha=0.8, \gamma=0.4424, n=0.1, \beta=0.138485, m=0.263075, \epsilon=0.01.$ Taking the initial value $(0.25, 0.1375),$ we observe that the trajectory spirals outward in a counterclockwise direction. On the other hand, when the initial value is chosen as $(0.25, 0.13),$ the trajectory spirals inward in a counterclockwise direction. According to  Poincaré-Bendixson theorem, there  exists an unstable limit cycle between these two trajectories. See Figure (\ref{unstable}).
\end{example}
\subsection{Cyclicity of slow-fast cycles}
Let $\Gamma(s)$ be slow-fast cycles without head, then it consists of the following segments:
\begin{equation*}
    \Gamma(s):=\left\{ (x,F(x))\mid x\in (\alpha_s, \omega_s) \right\} \cup \left\{(x,y_M-s)\mid x\in (\alpha_s, \omega_s), s\in (0,y_M-\hat{y}) \right\}.
\end{equation*}
In order to study the cyclicity  $Cycl(\Gamma(s))$ of the slow-fast cycles $\Gamma(s)$, we  define the slow divergence integral $I(s)$ (see \cite{Lizhu} and \cite{zzhu}). $I(s)$ is  defined as
\begin{equation}
\begin{split}
I(s,\lambda_0)&=\int_{\omega_s}^{\alpha_s} \frac{\partial f}{\partial x}(x, F(x), \lambda_0, 0)\frac{F'(x)}{g(x, F(x), \lambda_0, 0)}dx,  s\in(0, y_M-\hat{y}).\\
\end{split}
\label{sdi}
\end{equation}
 We deal with the  number of zeros of $I(s)$ to study the cyclicity  $Cycl(\Gamma(s))$.
Define $\sigma(x) \in [\alpha_s, \omega_s]$, such that for $x \in [x_M, \omega_s]$, $F(\sigma(x)) = F(x)$. Setting
\begin{equation}
    h(x)=\frac{\partial_x f(x,F(x),\lambda_0,0)}{g(x,F(x),\lambda_0,0)}.
\end{equation}
With direct calculations, we have
\begin{equation*}
    \begin{split}
    x&=\frac{1}{2}\left(1-m-n-y+\sqrt{(1-m-n-y)^2-4n(m+y)}\right),\\
    \sigma(x)&=\frac{1}{2}\left(1-m-n-y-\sqrt{(1-m-n-y)^2-4n(m+y)}\right),   
    \end{split}
\end{equation*}
then
\begin{equation*}
    I(s,\lambda_0)=\int_{y_M}^{y_M-s} h(\sigma(x))-h(x)\mid_{x=F^{-1}(y)}dy.
\end{equation*}
Set 
\begin{equation*}
 \Psi(x)=\frac{m-\sqrt{m}+x}{(m+x)^2\left(\alpha x-\beta-\gamma F(x)\right)\left(\frac{x}{m+x}-n-x\right)},   
\end{equation*}
then 
\begin{equation*}
     h(\sigma(x))-h(x)=\Psi(\sigma(x))\Psi(x)(\sigma(x)-x)m^{\frac{2}{3}}F(x)\Phi(y), 
\end{equation*}
where
\begin{equation*}
    \Phi(y)=(\alpha+\gamma)y+\gamma+n(\alpha+\gamma)-\sqrt{m}\alpha-m(\alpha+\gamma).
\end{equation*}
Recall that $\sigma(x)\in(\Bar{x}_1,x_M),x\in(x_M,\Bar{x}_2),$ 
then $\Psi(\sigma(x))\Psi(x)(\sigma(x)-x)<0.$ On the other hand,  
\begin{equation*}
\begin{split}
    \Phi(0)&=n(\alpha+\gamma)+(1-m)\gamma-\sqrt{m}\alpha(1+\sqrt{m}),\\
     \Phi(y_M)&=\alpha+2\gamma-(3\alpha+2\gamma)\sqrt{m},
\end{split}
 \end{equation*}
 and 
 \begin{equation*}
     \Phi(y)=0\Longleftrightarrow y=\frac{\sqrt{m}\alpha+(m-n)(\alpha+\gamma)-\gamma}{\alpha+\gamma}:=y_0.
 \end{equation*}
From
\begin{equation*}
  \Phi'(y)=\alpha+\gamma>0,
\end{equation*} it yields that
$\Phi(y)$ is a monotonically increasing function. Therefore, $\Phi(y_M-s) < \Phi(y) < \Phi(y_M)$
\begin{figure}[!ht]
\centering
\begin{minipage}[t]{0.3\textwidth}
\centering
\includegraphics[width=4cm]{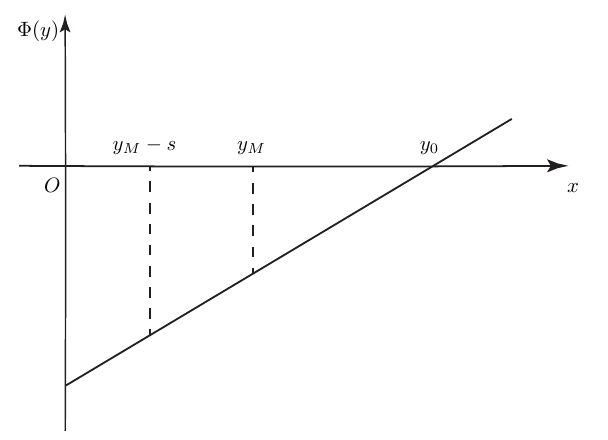}
\put(-75,-20){(a)}
\end{minipage}
\begin{minipage}[t]{0.3\textwidth}
\centering
\includegraphics[width=4cm]{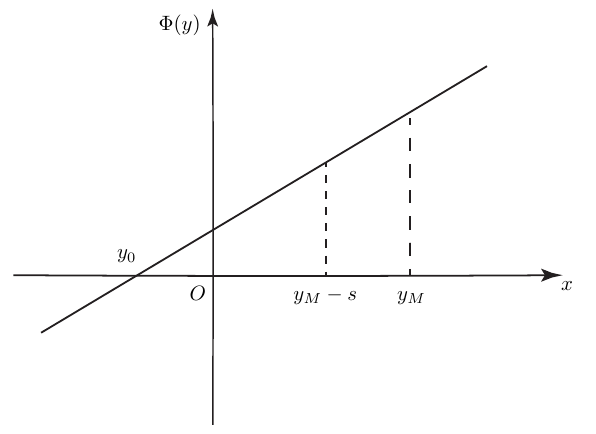}
\put(-60,-20){(b)}
\end{minipage}
\begin{minipage}[t]{0.3\textwidth}
\centering
\includegraphics[width=4cm]{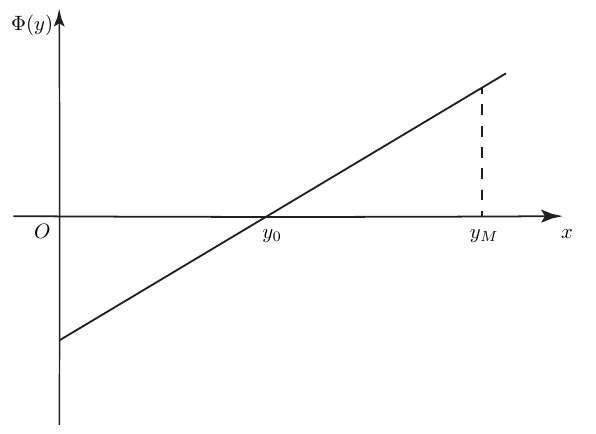}
\put(-65,-20){(c)}
\end{minipage}
\caption{
The relative positions of $y_0$ and $y_M.$}
\label{YY}
\end{figure}
and 
(1) if $m\ge\left(\frac{\alpha+2 \gamma}{3 \alpha+2 \gamma}\right)^2,$ i.e. $y_0\ge y_M,$  then  $\Phi(y)<0.$ See Fig.(\ref{YY})(a).

(2) if  $m<\left(\frac{\alpha+2 \gamma}{3 \alpha+2 \gamma}\right)^2,$ i.e. $y_0 < y_M,$  See Fig. (\ref{YY})(b)(c).   In this case,  to determine the sign of $\Phi(y)$,   we need to consider the relative position of $y_0$ with respect to 0. Then
 
a) if $y_0 = 0$, then $y_m > 0$, and when $y \in (y_M - s, y_M)$, $\Phi(y) > 0$.

b) if $y_0 < 0$, then $\Phi(y) > 0$, see Fig. (\ref{YY})(b).

c) if $y_0 > 0$, then by the expression for $y_0$ and  $0 < n < (1-\sqrt{m})^2$, we have $m < \left(\frac{\alpha+2 \gamma}{3 \alpha+2 \gamma}\right)^2$. This contradicts $m > \left(\frac{\alpha+2 \gamma}{3 \alpha+2 \gamma}\right)^2,$ see Fig. (\ref{YY})(c).

Therefore, when $m<\left(\frac{\alpha+2 \gamma}{3 \alpha+2 \gamma}\right)^2$, $y_0\le 0$, $\Psi(y)>0$, and when $m\ge\left(\frac{\alpha+2 \gamma}{3 \alpha+2 \gamma}\right)^2$, $\Psi(y)<0.$
\begin{theorem}
For system (\ref{2}) with  $\gamma=\frac{\beta +\alpha  m-\alpha  \sqrt{m}}{-m+2 \sqrt{m}+n-1}$, $\bigtriangleup_1>0,$ and $1-m-n>0$, $Cycl(\Gamma(s)) \le 1$. 
\end{theorem}

\section{Discussion}
In this paper, we extend the Lyapunov coefficient formula for singular Hopf bifurcations from \cite{krupa2} for general geometric singular perturbation systems. This extension addresses the case when the $O(1)$-order Lyapunov coefficient $A=0$, a situation not covered in the original work, see \cite{krupa2}. We also discuss the singular Hopf bifurcation  when $A=0$. To illustrate its practical application, we investigate the dynamical behavior of a predator-prey model with Allee effects. Assume that the prey's reproductive rate is significantly higher than that of the predator, we reduce this ecological system to a slow-fast system with a small parameter. In the multi-scale framework, we obtain the  cyclicity  of singular Hopf bifurcation and canard explosion bifurcation curves. Additionally, utilizing the proposed method for calculating the first Lyapunov coefficient, we derive that when $A=0$, the system exhibits a subcritical singular Hopf bifurcation, resulting in an unstable limit cycle. Numerical examples confirm  our results. Finally, using the slow divergence integral theory, we demonstrate that the cyclicity of slow-fast limit periodic set  is 1.

It is worth noting that our extended method can be applied to find higher-order approximations or higher-order Lyapunov coefficient formulas. Specifically, if the first Lyapunov coefficient $L_{\epsilon}\equiv 0$, indicating the occurrence of a degenerate Hopf bifurcation, what is the second-order Lyapunov coefficient?  Furthermore, in the derived first-order Lyapunov coefficient formula, $\omega_1=0$ is equivalent to $A=0$ in \cite{krupa2}. Under this condition, if we assume $\omega_2|{\omega_1=0}<0$ or $\omega_2|{\omega_1=0}>0$, will the canard explosion still occur? What is the mechanism behind it? This is another crucial question that we plan to explore in our future research.

\section*{Data availability}
No data was used for the research described in the article.

 \bibliographystyle{elsarticle-num} 
 \bibliography{cas-refs}





\end{document}